\documentclass{gtart}

\input gtmonout
\volumenumber{2}
\volumeyear{1999}
\volumename{Proceedings of the Kirbyfest}
\pagenumbers{215}{232}
\papernumber{12}
\received{20 January 1999}\revised{28 July 1999}
\published{18 November 1999}

\newtheorem*{lem}{Lemma}
\newtheorem*{thm}{Theorem}
\newtheorem*{tthm}{Trace Theorem}
\theoremstyle{remark}

\newtheorem*{rem}{Remark}

\begin{document}

\title{Right integrals and invariants of three--manifolds}
\asciititle{Right integrals and invariants of three-manifolds}

\author{Louis H Kauffman}

\address{Department of Mathematics,
Statistics and Computer Science\\University of Illinois at
Chicago\\851 South Morgan Street, Chicago, IL 60607-7045, USA}
\email{kauffman@uic.edu}

\begin{abstract}
This paper gives a summary of our approach to invariants
of three manifolds via right integrals on finite dimensional
Hopf algebras and their relation to the Kirby calculus. 

{\em It gives the author great pleasure to dedicate this paper to
Rob Kirby on his sixtieth birthday.}
\end{abstract}
\asciiabstract{This paper gives a summary of our approach to
invariants of three manifolds via right integrals on finite
dimensional Hopf algebras and their relation to the Kirby calculus.}


\makeshorttitle

\section{Introduction}

This paper is a capsule summary of our approach to invariants
of three manifolds via right integrals on finite dimensional
Hopf algebras. The aim of this paper is to quickly show the
reader how natural it is to consider right integrals in this
context. The defining property of a right integral is a
categorical algebra translation of the handle-sliding move in
the Kirby calculus.

We show  in \cite{KandR} that invariants defined in terms of
right integrals, as considered in this paper,  are distinct
from the invariants of Reshetikhin and Turaev \cite{RT}.  We
show  that the Hennings invariant (defined via these
integrals)  is non-trivial for the quantum group
$U_{q}(sl_{2})'$ when  $q$  is an fourth root of unity.  The
Reshetikhin Turaev invariant is trivial at this quantum group
and root of unity.  The Hennings invariant  distinguishes all
the Lens spaces  $L(n,1)$  from one another at this root of
unity. This proves that there is non-trivial topological
information in the non-semisimplicity  of
$U_{q}(sl_{2})'$. This non-triviality  result has also been
obtained by Ohtsuki
\cite{Ohtsuki}. 

The reader interested in comparing the approach of this paper
with other ways to look at quantum link invariants will enjoy
looking at  the references
\cite{KandP}, \cite{Kerler}, \cite{KirbyMelvin},
\cite{Lawrence}, \cite{Luba},
\cite{Radford-trace}, \cite{Radford-gen}, \cite{Resh},
\cite{RT},
\cite{Witten}.
 In particular, the method we use to write link invariants
directly in relation to a Hopf algebra is an analog of the
construction in \cite{Lawrence} and it is a generalization of
the formalism of \cite{Resh} and \cite{RT}.  The papers
\cite{Kerler}, \cite{Luba}, \cite{Kuperberg1},
\cite{Kuperberg2} consider categorical frameworks that also use
right integrals on Hopf algebras.

The paper is organized as follows.  Section 2 recalls the Hopf
algebra background. Section 3  recalls the category associated
with a Hopf algebra. Section 4 weaves in the Tangle category
and constructs the functor from the Tangle category to the Hopf
algebra category. Section 5 describes the invariants of three
manifolds. 

\noindent {\bf Acknowledgement}\qua    I thank the National
Science Foundation for support of this research under NSF Grant
DMS-9205277 and the NSA for  partial support under grant number
MSPF-96G-179. 

\section{Algebra} Recall that a Hopf algebra  $A$
\cite{Sweedler}   is a bialgebra over a commutative ring  $k$
that has an associative multiplication and  a coassociative
comultiplication and  is equipped with a counit, a unit and an
antipode. The ring $k$ is usually taken to be a field. $A$ is
an algebra with multiplication $m\co A \otimes A \longrightarrow
A$. The associative law for m is expressed by the equation $m(m
\otimes1_{A}) = m(1_{A} \otimes m)$ where
$1_{A}$ denotes the identity map on A. 

The coproduct $\Delta \co A \longrightarrow A \otimes A$  is an
algebra homomorphism and is coassociative in the sense that
$(\Delta \otimes 1_{A})\Delta = (1_{A}
\otimes \Delta) \Delta.$ 

The unit is a mapping from  $k$  to $A$  taking  $1$ in $k$ to
$1$ in $A$, and thereby defining an action of $k$ on $A.$  It
will be convenient to just identify the units in $k$ and in
$A$, and to ignore the name of the map that gives the unit.

The counit is an algebra mapping from  $A$  to  $k$ denoted by
$\epsilon \co A
\longrightarrow  k.$ The following formulas for the counit
dualize the structure inherent in the unit: $(\epsilon \otimes
1_{A}) \Delta = 1_{A} = (1_{A} \otimes
\epsilon) \Delta.$ 

It is convenient to write formally 
$$\Delta (x) = \sum x_{1}
\otimes x_{2}
\in A
\otimes A$$ 
to indicate the decomposition of the coproduct of
$x$ into a sum of first and second factors in the two--fold
tensor product of $A$ with itself.  We shall often drop the
summation sign and write 
$$\Delta (x) =  x_{1} \otimes x_{2}.$$
The antipode is a mapping  $s\co A \longrightarrow A$ satisfying
the equations
$m(1_{A} \otimes s) \Delta (x) = \epsilon (x)1,$  and $m(s
\otimes 1_{A})
\Delta (x)= \epsilon (x)1$ where 1 on the right hand side of
these equations denotes the unit of $k$ as identified with the
unit of $A.$ It is a consequence of this definition  that
$s(xy) = s(y)s(x)$  for all  $x$ and  $y$  in  A. 

A  quasitriangular Hopf algebra  $A$   \cite{Drinfeld}  is a
Hopf algebra with an element  $\rho \in A \otimes A$ satisfying
the following equations:

1)  $\rho \Delta = \Delta' \rho$   where  $\Delta'$
is the composition of  $\Delta$  with the map on $A \otimes A$
that switches the two factors.

\noindent 2) $$\rho_{13} \rho_{12} = (1_{A} \otimes \Delta)
\rho,$$ $$\rho_{13}
\rho_{23} = (\Delta \otimes 1_{A})\rho.$$ 

\begin{rem}The symbol $\rho_{ij}$ denotes the
placement of the first and second tensor factors of $\rho$  in
the $i$ and $j$ places in a triple tensor product.  For
example, if $\rho = \sum e \otimes e'$ then  $$\rho_{13} = \sum
e
\otimes 1_{A} \otimes e'.$$ 
These conditions imply   that  $\rho$   has an inverse, and that
$$ \rho^{-1} = (1_{A} \otimes s^{-1}) \rho = (s \otimes 1_{A})
\rho .$$
It follows easily from the axioms of the quasitriangular Hopf
algebra that
$\rho$  satisfies the Yang--Baxter equation
$$\rho_{12} \rho_{13} \rho_{23} = \rho_{23} \rho_{13}
\rho_{12}.$$\end{rem}

A less obvious fact about quasitriangular Hopf algebras is that
there exists an element  $u$  such that $u$  is invertible and
$s^{2}(x) = uxu^{-1}$  for all
$x$ in  $A.$ In fact, we may take   $u = \sum s(e')e$  where
$\rho = \sum e
\otimes e'$.  This result is due to Drinfeld \cite{Drinfeld}.   

An  element  $G$  in a Hopf algebra is said to be  {\em
grouplike}   if
$\Delta (G) = G \otimes G$  and  $\epsilon (G)=1$  (from which
it follows that
$G$ is invertible and $s(G) = G^{-1}$).    A quasitriangular
Hopf algebra is said to be a {\em ribbon Hopf algebra}
\cite{RTG}, \cite{KR} if there exists a grouplike element $G$
such that  (with  $u$ as in the previous paragraph) $v =
G^{-1}u$    is in the center of $A$  and  $s(u) =
G^{-1}uG^{-1}$.   We call G a {\em special  grouplike} element of A.

Since  $v=G^{-1}u$  is central,  $vx=xv$ for all x  in A.
Therefore
$G^{-1}ux = xG^{-1}u.$ We know that $s^{2}(x) = uxu^{-1}.$
Thus  $s^{2}(x) =GxG^{-1}$  for all $x$  in $A.$ Similarly,
$s(v) = s(G^{-1}u) = s(u)s(G^{-1})=G^{-1}uG^{-1}G =G^{-1}u=v.$
Thus the square of the antipode is represented as conjugation by
the special grouplike element in a ribbon Hopf algebra, and the
central element
$v=G^{-1}u$  is invariant under the antipode.

\section{Categories and Functors} We will recall in capsule
form the method for obtaining invariants of knots, links and
3--manifolds that we used in
\cite{KandR}, \cite{K-Hopf} and \cite{KRS}.  This method
produces a functor
$$F\co Tang \longrightarrow Cat(A)$$ where $Tang$ is the category
of unoriented tangles taken up to regular isotopy, and $Cat(A)$
is a natural category associated with a quasi-triangular Hopf
algebra.  We first recall the structure of $Cat(A)$, then show
how right integrals on the Hopf algebra lead to invariants of
three manifolds. 

In this recollection, we start by describing the association of
a category with an algebra $A$, not necessarily a Hopf algebra,
and then add structure to the category corresponding to
properties of a Hopf algebra. By the end of the discussion, the
category will have been built. 

At the outset $A$ is an algebra with a commutative ground ring
$R$, so that $A$ is a module over $R.$ $Cat(A)$ is an
(associative tensor) category whose objects are generated under
tensor product by two objects $k$ and $V$ such that $k
\otimes k = k$ and $k \otimes V = V \otimes k = V,$ and whose
morphisms are generated using composition and tensor product,
from basic morphisms that will be defined in this section.  The
object $V$ is strictly abstract. The object $k$ will be
identified with the ground ring $R$ as described below.
Functors on this category will often associate to $V$ a
representation space of the algebra $A.$

To each  element $a$ of the algebra $A$ there will be a
generating morphism (still denoted by $a$)  $$a\co V
\longrightarrow V.$$  Multiplication in the algebra corresponds
to composition of the corresponding morphisms.  The identity
morphism
$1_{V}$ corresponds to the identity element $1$  in the algebra
$A$   In the same way, elements of the tensor powers
$A^{\otimes n}$ of the algebra correspond to morphisms with
domain and range $V^{\otimes n}.$ 

Morphisms in $Cat(A)$ corresponding to elements of the algebra
are denoted by vertical upward pointing arrows (with the domain
below the range).  We shall customarily denote the presence of
the algebra element $a$ by a ``bead" (a node) on the arrow with
the label $a.$ An unlabelled arrow without a bead corresponds
to the identity morphism $1_{V}.$ See Figure 1.

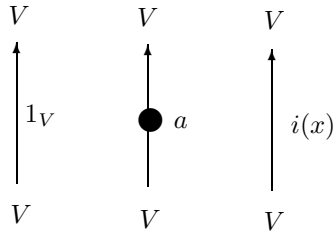
\begin{figure}[ht!]\vglue -3mm\cl{%
\setlength{\unitlength}{0.92pt}
\begin{picture}(160,124) \thinlines\small
\put(73,61){\circle*{10}} \put(130,49){\makebox(20,19){$i(x)$}}
\put(21,55){\makebox(14,16){$1_{V}$}}
\put(78,52){\makebox(15,16){$a$}}
\put(115,93){\makebox(18,19){$V$}}
\put(115,10){\makebox(18,19){$V$}}
\put(63,94){\makebox(18,19){$V$}}
\put(64,11){\makebox(18,19){$V$}}
\put(10,95){\makebox(18,19){$V$}}
\put(11,13){\makebox(18,19){$V$}}
\put(123,32){\vector(0,1){58}} \put(72,34){\vector(0,1){58}}
\put(18,35){\vector(0,1){58}} \end{picture}}\vglue -5mm
\caption{Labelled Arrow}
\end{figure}

The algebra $A$ has a base ring $R$ (usually a field).  The
object $k$ in the category is identified with the base ring
$R$.  If $i\co  k \longrightarrow A$ denotes the unit in the
algebra (ie, $i(x) = x1$ where $x$ is in $k$ and
$1$ is the unit in the algebra) then elements of $i(k)$ form a
special set of morphisms in the category since $$i(x) a \otimes
b = a \otimes i(x) b$$ for any $x$ in $k$ and $a$ and $b$ in
$A$.  Accordingly, we shall denote a morphism in
$Cat(A)$ that corresponds to an element $w$ in $i(k)$ by an
arrow {\em without a bead} but labelled with $w$. In this way
we can move the label $w$ freely and when it leaves a given
arrow, that arrow (if otherwise unlabelled) reverts to an
identity arrow as it should.  See Figure 1. 

We are free to utilize the internal structure of the base ring
$k$ when that is appropriate. It is a curious feature of this
category that one of its objects is also the holder of an
important subcategory (the morphisms generated by
$i(k)$). In other respects the base object $k$ can be regarded
in a purely external categorical way. 

Every Hopf algebra has an antipode $s\co  A \longrightarrow A$
with the defining property that
$$m(s \otimes 1_{A})(\Delta (x)) = m(s \otimes 1_{A})(\Delta
(x)) = i \epsilon (x)$$
\noindent where $m\co  A \otimes A \longrightarrow A$ is the
multiplication in the algebra,  $\epsilon \co  A \longrightarrow
k$ is the counit,  and $\Delta\co  A
\longrightarrow A \otimes A$ is the coproduct.  A key property
of the antipode is that it is an antimorphism of the algebra to
itself. It reverses the order of multiplication so that
$$s(ab) = s(b)s(a)$$
\noindent for all $a$ and $b$ in $A.$

The structure of the antipode is reflected  in the category
$Cat(A)$ by new morphisms that we now describe as
$$Cup\co  k \longrightarrow V \otimes V$$
\noindent and
$$Cap\co   V \otimes V \longrightarrow k.$$
\noindent The key property of these morphisms is their
behaviour with respect to the antipode:
$$ (a \otimes 1_{V}) \circ Cup =  (1_{V} \otimes s(a)) \circ
Cup, $$
$$ Cap \circ (1_{V} \otimes a) =  Cap \circ ( s(a) \otimes
1_{V}).$$
\noindent These morphisms and all the remaining morphisms that
we construct in this category are represented by (immersed)
curves in the plane. See Figure 2 for the diagrammatic
depiction of $Cup$, $Cap$ and these ``bead sliding" properties.
It is convenient to think of these identities as permissions to
slide the beads past maxima or minima in the diagrams for the
morphisms. Upon sliding across a maximum in a counterclockwise
direction, a bead has the antipode applied to it.

\begin{figure}[ht!]
\cl{\small    \setlength{\unitlength}{0.92pt}
\begin{picture}(264,142) \thinlines
\put(120,22){\makebox(16,15){$=$}}
\put(116,99){\makebox(16,15){$=$}}
\put(236,23){\makebox(18,14){$s(a)$}}
\put(10,21){\makebox(15,14){$a$}}
\put(180,107){\makebox(18,14){$s(a)$}}
\put(87,108){\makebox(15,14){$a$}}
\put(224,29){\circle*{10}} \put(32,28){\circle*{10}}
\put(167,115){\circle*{10}}
\put(81,115){\circle*{10}} \put(59,39){\oval(56,58)[b]}
\put(197,40){\oval(56,58)[b]} \put(194,103){\oval(56,58)[t]}
\put(55,103){\oval(56,58)[t]} \end{picture}}
\caption{Bead Sliding for Cups and Caps}
\end{figure}
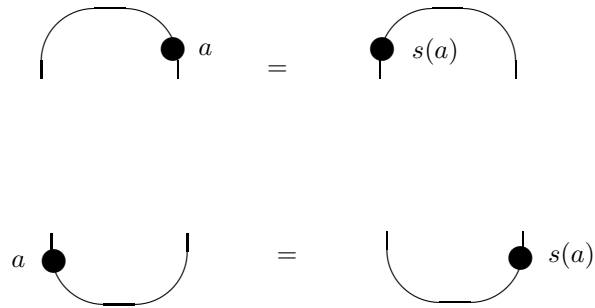

The properties of $Cup$ and $Cap$ are augmented with the {\em
topological} identities
$$(Cup \otimes 1_{V})(1_{V} \otimes Cap) = 1_{V},$$
$$(1_{V} \otimes Cup)(Cap \otimes 1_{V}) =1_{V}.$$
\noindent See Figure 3.  These identities assert that $Cup$ and
$Cap$ behave just as topological maxima and minima: when
composed without any intervening morphisms they cancel to the
identity.

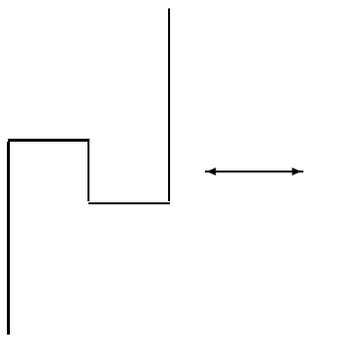
\begin{figure}[ht!]\vglue-2mm
\cl{\small    \setlength{\unitlength}{0.92pt}
\begin{picture}(161,158) \thinlines
\put(122,79){\vector(-1,0){29}} \put(100,79){\vector(1,0){33}}
\put(149,148){\line(0,-1){138}} \put(12,91){\line(0,-1){79}}
\put(78,146){\line(0,-1){79}} \put(12,92){\line(1,0){33}}
\put(45,66){\line(1,0){33}} \put(45,92){\line(0,-1){25}}
\end{picture}}\vglue-2mm
\caption{Cup Cap Cancellation}
\end{figure}

It follows from the antipodal and topological behaviour of
$Cup$ and $Cap$ that the morphism corresponding to $s(a)$ is
represented by composing on the top and the bottom of $a$ with
a $Cup$ and a $Cap$, as shown in Figure 4.  A proof by bead
sliding is shown in that Figure. 

\begin{figure}[ht!]
\cl{\small    \setlength{\unitlength}{0.92pt}
\begin{picture}(302,172) \thinlines
\put(55,74){\makebox(17,21){$=$}}
\put(194,74){\makebox(17,21){$=$}}
\put(266,80){\makebox(14,16){$a$}}
\put(89,75){\makebox(16,19){$s(a)$}}
\put(10,76){\makebox(16,19){$s(a)$}} \put(258,88){\circle*{10}}
\put(118,84){\circle*{10}} \put(290,159){\line(0,-1){79}}
\put(224,91){\line(0,-1){79}} \put(274,87){\oval(34,42)[b]}
\put(241,87){\oval(34,42)[t]} \put(37,85){\circle*{10}}
\put(134,85){\oval(34,42)[t]} \put(167,85){\oval(34,42)[b]}
\put(117,89){\line(0,-1){79}} \put(183,157){\line(0,-1){79}}
\put(36,162){\line(0,-1){152}} \end{picture}}
\caption{Bead Sliding for the Antipode}
\end{figure}
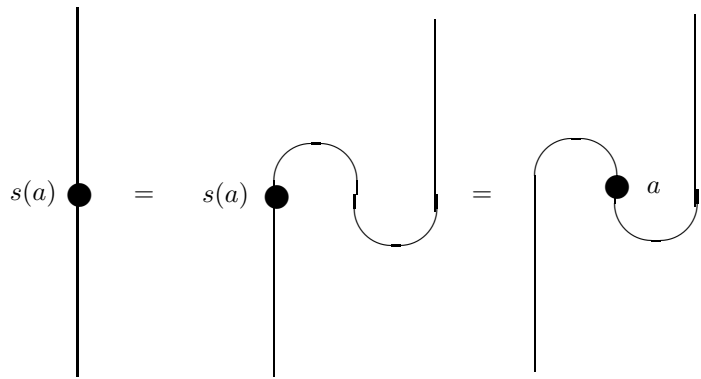 

Formally, we have the result:

\begin{lem}Let $a$ be an element of $A$ viewed as
a morphism from
$V$ to $V$.  Then  $s(a) = (1_{V} \otimes Cup)(1_{V} \otimes a
\otimes 1_{V})(Cap
\otimes 1_{V})$.\end{lem}

\begin{proof}\quad $(1_{V} \otimes Cup)(1_{V} \otimes a
\otimes 1_{V})(Cap
\otimes 1_{V})$
\nl\hbox{}\qquad\quad$=(1_{V} \otimes Cup)(1_{V} \otimes 1_{V}
\otimes s(a))(Cap
\otimes 1_{V})$
\nl\hbox{}\qquad\quad$ =(1_{V} \otimes Cup) (Cap \otimes 1_{V})s(a)$
\nl\hbox{}\qquad\quad$=1_{V}s(a) = s(a).$ \end{proof}

In fact, the cup and cap morphisms provide an interpretation of
$s(f)$ for any morphism $f\co  V \longrightarrow V$ in the
category by the replacement of $a$ by
$f$ in the diagram in Figure 4, and in the formula of this
lemma. Thus
$$s(f) = (1_{V} \otimes Cup)(1_{V} \otimes f \otimes 1_{V})(Cap
\otimes 1_{V}).$$

One last type of morphism completes
that catalog for $Cat(A).$ This morphism is a formal
permutation
$$P\co  V \otimes V \longrightarrow V \otimes V$$
\noindent with the property
$$P \circ (a \otimes b) = (b\otimes a) \circ P.$$
\noindent We represent $P$ by a pair of crossed arrows as shown
in Figure 5.

\begin{figure}[ht!]
\cl{\small    \setlength{\unitlength}{0.92pt}
\begin{picture}(245,174) \thinlines
\put(83,121){\makebox(72,20){Permutation}}
\put(85,21){\makebox(19,17){$b$}}
\put(12,23){\makebox(16,17){$a$}}
\put(139,55){\makebox(19,17){$b$}}
\put(219,54){\makebox(16,17){$a$}}
\put(161,120){\makebox(24,26){$P$}}
\put(36,29){\circle*{10}} \put(75,29){\circle*{10}}
\put(168,64){\circle*{10}}
\put(209,63){\circle*{10}} \put(11,162){\line(4,-3){85}}
\put(11,98){\line(4,3){85}} \put(12,76){\line(4,-3){85}}
\put(12,12){\line(4,3){85}} \put(146,17){\line(4,3){85}}
\put(146,81){\line(4,-3){85}} \end{picture}}
\caption{The Permutation }
\end{figure}
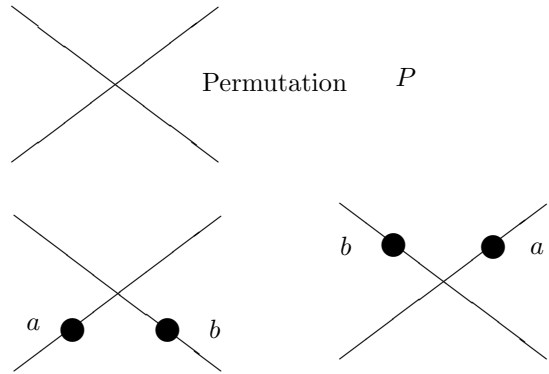

It is assumed that $P$ satisfies the usual axioms for
permutation plus a compatibility condition with respect to the
$Cup$ and the $Cap$ as shown in Figure 6.

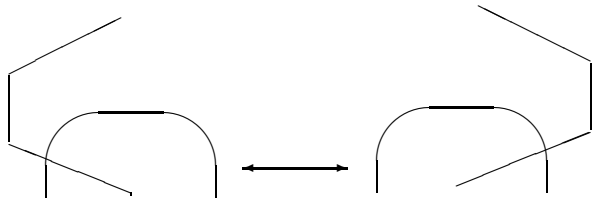
\begin{figure}[ht!]\vglue-2mm
\cl{\small    \setlength{\unitlength}{0.92pt}
\begin{picture}(263,136) \thinlines
\put(139,57){\vector(-1,0){31}} \put(112,57){\vector(1,0){39}}
\put(251,101){\line(-2,1){46}} \put(251,73){\line(0,1){27}}
\put(196,50){\line(5,2){55}} \put(12,96){\line(2,1){46}}
\put(12,68){\line(0,1){27}} \put(62,47){\line(-5,2){50}}
\put(62,46){\line(0,1){1}} \put(198,47){\oval(70,70)[t]}
\put(62,45){\oval(70,70)[t]} \end{picture}}\vglue-8mm
\caption{Permutation Switchback}\end{figure}

The Hopf algebra $A$ is endowed with a coproduct
$$\Delta\co  A  \longrightarrow A \otimes A.$$
\noindent This map is indicated algebraically via the formula
$$\Delta(a)  = \Sigma a_{1} \otimes a_{2}$$
\noindent where an implicit summation is sometimes indicated
without the presence of the summation symbol.  In diagrams the
line with a bead labelled $a$ on it is replaced by a two
parallel lines (for the tensor product) with the left line's
bead labelled $a_{1}$ and the right line's bead labelled
$a_{2}.$  See Figure 7.

\begin{figure}[ht!]
\cl{\small    \setlength{\unitlength}{0.92pt}
\begin{picture}(174,124) \thinlines
\put(72,53){\makebox(18,15){$=$}}
\put(146,51){\makebox(18,17){$a_{2}$}}
\put(111,51){\makebox(19,16){$a_{1}$}}
\put(61,53){\makebox(13,18){$)$}}
\put(10,50){\makebox(29,24){$\Delta($}}
\put(104,60){\circle*{10}}
\put(94,93){\makebox(18,19){$V$}}
\put(129,10){\makebox(18,19){$V$}}
\put(103,33){\vector(0,1){58}} \put(139,60){\circle*{10}}
\put(129,92){\makebox(18,19){$V$}}
\put(96,11){\makebox(18,19){$V$}}
\put(138,33){\vector(0,1){58}} \put(42,35){\vector(0,1){58}}
\put(34,12){\makebox(18,19){$V$}}
\put(33,95){\makebox(18,19){$V$}}
\put(47,53){\makebox(15,16){$a$}} \put(43,62){\circle*{10}}
\end{picture}}\vglue-4mm
\caption{Diagramming the coproduct}\end{figure}
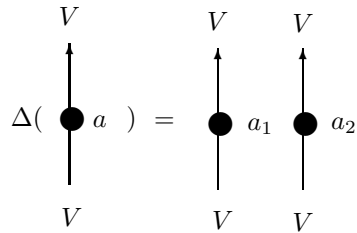

Clearly $\Delta$ must simply double all diagrammatic lines to
be defined on the structural morphisms in $Cat(A)$.  This is
illustrated in Figure 8.

\begin{figure}[ht!]
\cl{\small    \setlength{\unitlength}{0.92pt}
\begin{picture}(271,189) \thinlines
\put(203,34){\makebox(58,30){$=G^{-1}$}}
\put(51,75){\line(1,0){36}}
\put(142,138){\oval(32,32)} \put(141,138){\oval(16,18)}
\put(204,141){\oval(32,32)} \put(243,142){\oval(32,32)}
\put(161,129){\makebox(15,15){$=$}}
\put(188,178){\line(0,-1){75}}
\put(132,177){\line(0,-1){75}} \put(125,177){\line(0,-1){75}}
\put(227,179){\line(0,-1){75}} \put(103,50){\oval(32,32)}
\put(16,176){\line(0,-1){75}} \put(33,138){\oval(32,32)}
\put(51,74){\line(0,-1){64}} \put(87,20){\line(1,0){36}}
\put(190,87){\line(0,-1){75}} \put(175,50){\oval(32,32)}
\put(23,126){\makebox(21,23){$G$}}
\put(10,38){\makebox(28,22){$s(G)=$}}
\put(135,39){\makebox(20,21){$=$}}
\put(79,125){\makebox(29,24){$\Delta(G) =$}}
\put(87,74){\line(0,-1){53}} \put(124,20){\line(0,1){68}}
\end{picture}}
\caption{Delta on Diagrams}
\end{figure}
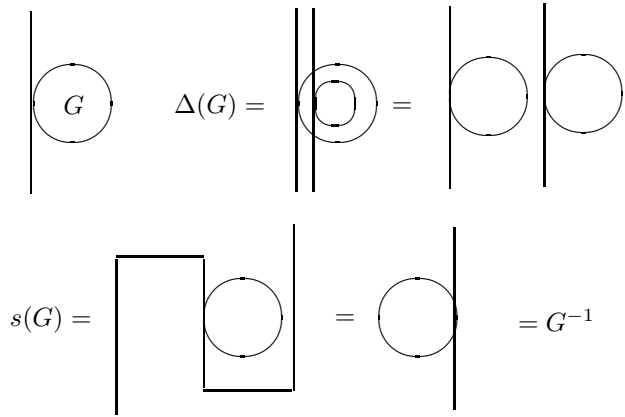

\subsection{Curls} Before going to the tangle category it is
useful to delineate special morphisms in $Cat(A)$ that we call
``curls". 

A {\em curl} is a morphism of the form $G=(1_{V} \otimes
Cup)P(1_{V} \otimes Cap)$ or $G^{-1} = (Cup \otimes 1_{V})P(Cap
\otimes 1_{V})$.  The second curl is denoted by $G^{-1}$
because the equations $GG^{- 1} = 1_{V} = G^{-1}G$ hold via a
regular homotopy of the corresponding plane curves. (The
regular homotopy is implicitly shown in Figure 9.) These
equations are the categorical analog of the so-called ``Whitney
trick".  The Whitney Graustein Theorem \cite{Whitney} tells us
that any morphism from $V$ to $V$ with no elements of $A$ in its
factorization is equivalent to an integer power of $G$.

The element $G$ is a formal grouplike element in the sense that
$\Delta(G) = G
\otimes G$ as illustrated in Figure 8. 

In terms of the morphism $G$ we can represent the square of the
antipode in
$Cat(A)$ by the formula  $$s^{2}(a) = GaG^{-1}$$ where this
formula denotes a composition of morphisms in $Cat(A).$  In
many cases the element $G$ can be identified with a specific
grouplike element in the Hopf algebra $A$. The category takes
the lead here in expressing a preference for such a formula. See
Figure 9.

\begin{figure}[ht!]\vglue-2mm
\cl{\small    \setlength{\unitlength}{0.92pt}
\begin{picture}(328,168) \thinlines
\put(28,84){\makebox(17,19){$a$}}
\put(214,87){\makebox(22,19){$=$}}
\put(103,87){\makebox(22,19){$=$}}
\put(256,10){\makebox(52,28){$GaG^{-1}$}}
\put(28,10){\makebox(35,29){$s^{2}(a)$}}
\put(302,76){\line(0,1){12}}
\put(257,96){\line(0,1){15}} \put(302,89){\line(-3,-1){59}}
\put(316,116){\line(-3,-1){59}} \put(60,85){\line(0,0){0}}
\put(165,74){\line(0,0){0}} \put(302,75){\line(-1,0){22}}
\put(279,111){\line(-1,0){22}} \put(279,110){\line(0,-1){35}}
\put(280,92){\circle*{10}} \put(316,116){\line(0,1){40}}
\put(242,29){\line(0,1){40}} \put(131,29){\line(0,1){40}}
\put(205,116){\line(0,1){40}} \put(131,68){\line(1,1){59}}
\put(146,57){\line(1,1){59}} \put(169,92){\circle*{10}}
\put(168,110){\line(0,-1){35}} \put(168,111){\line(-1,0){22}}
\put(191,75){\line(-1,0){22}} \put(146,111){\line(0,-1){53}}
\put(191,128){\line(0,-1){53}} \put(71,76){\line(0,0){0}}
\put(86,158){\line(0,-1){99}} \put(-34,87){\line(0,0){0}}
\put(12,129){\line(0,-1){99}} \put(13,129){\line(1,0){59}}
\put(27,58){\line(1,0){59}} \put(72,129){\line(0,-1){53}}
\put(27,112){\line(0,-1){53}} \put(72,76){\line(-1,0){22}}
\put(49,112){\line(-1,0){22}} \put(49,111){\line(0,-1){35}}
\put(50,93){\circle*{10}} \end{picture}}\vglue -5mm
\caption{The Square of the Antipode}
\end{figure}
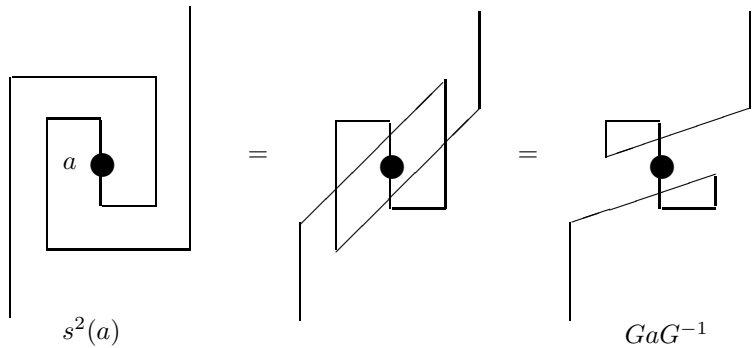 

It is now easy to see the following:

\begin{lem}Every morphism $x$ in $Cat(A)$ with
source and target both
$V$ can be written in the form $aG^{n}$ where $a \in{}A$ and
$n$ is an integer.\end{lem}

\begin{proof} This result is accomplished by sliding
all the algebra on the immersed curve that represents the
morphism to the bottom of the curve. The element $a$ is the
algebraic element so concentrated at the top of the curve. The
curve itself is regularly homotopic (fixing the ends at top and
bottom) to
$G^{n}$ for a unique integer $n$ (the Whitney degree of that
curve). Since regular homotopy is produced by the topological
axioms of the category, this shows that $x = aG^{n}$ as
desired.\end{proof} 

\begin{rem}Using the method of this proof, we
define $w(x) = a$ and
$d(x) = n$ so that $x = w(x)G^{d(x)}.$ \end{rem}

\section{The Tangle Category}

For the detailed definition of $Tang$ we refer the reader to
\cite{KRS}.  The objects will again be formally generated by
two (formal) symbols $k$ and $V$ as for $Cat(A).$  The
morphisms of  $Tang$ are generated by $Cup$, $Cap$, the
identity line (and its tensor products) and the {\em right and
left crossings}, maps of $V \otimes V$ to itself.  The
crossings are denoted  $R$ and $L$ as shown in Figure 10. $R$
denotes a crossing of two arcs so that the overcrossing line
goes upward from right to left. $L$ denotes a crossing of two
arcs so that the overcrossing line goes upward from left to
right. 

All diagrams of morphisms in $Tang$ are drawn in the plane with
no intersections of lines other than the formalized crossings.
The crossings are intended to be interpreted in sense of the
knot theorist. That is, one can {\em interpret} such a diagram
as representing an embedding in three--dimensional space of a
collection of curves. The crossing indicates that one curve
passes over another from the point of view of a planar
projection. 

Note that the permutation $P$ of the category $Cat(A)$ is
absent from $Tang.$ The crossings can be regarded as
generalizations of the permutation. We will shortly use this
idea in making a functor $F$ from $Tang$ to $Cat(A).$

\begin{figure}[ht!]\vglue -2mm
\cl{\small    \setlength{\unitlength}{0.92pt}
\begin{picture}(228,121) \thicklines
\put(159,10){\makebox(29,26){$L$}}
\put(44,10){\makebox(26,28){$R$}}
\put(178,68){\line(4,-3){36}} \put(167,77){\line(-4,3){38}}
\put(132,43){\line(4,3){85}} \put(59,79){\line(4,3){38}}
\put(49,71){\line(-4,-3){36}} \put(11,109){\line(4,-3){86}}
\end{picture}}\vglue -4mm
\caption{Right and Left Crossings}
\end{figure}
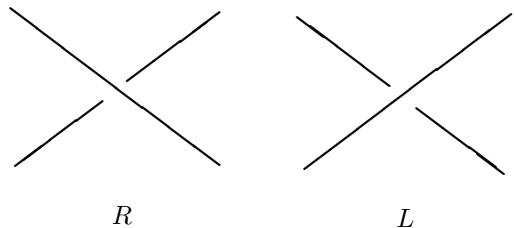

Every morphism in the tangle category is a composition of the
elementary morphisms $Cup$,$Cap$, $R$ and $L.$ The objects in
$Tang$ are identical to the objects in $Cat(A).$

While the objects in the tangle category are very simple, the
morphisms are quite complex. Each morphism in the tangle
category consists in a link diagram with free ends which is
transverse with respect to a given direction in the plane.
(This special direction will be called the {\em vertical}
direction.) The transversality of the diagram to this vertical
direction means that any given line perpendicular to the
vertical direction intersects the diagram either tangentially
at a maximum or a minimum, at non-zero angle for any other
strand. We shall further assume that any given perpendicular
intersects the diagram at most at one crossing.  With these
stipulations the free ends of the diagram occur at either its
top or its bottom (top and bottom taken with respect to the
designated vertical direction).We shall assume that all the top
ends occur along the same perpendicular, and that all the
bottom ends occur along another perpendicular to the vertical.
To each of these two rows of diagram ends is assigned a tensor
product of copies of $V$, one for each end.  In the tangle
category, the diagram is  a morphism from the lower tensor
product to the upper tensor product.   If the top of a diagram
has no free ends, then its range is $k$. If the bottom of a
diagram has no free ends, then its domain is $k$. A diagram is
said to be {\em closed} if it has no free ends. Thus a closed
diagram is a morphism from $k$ to
$k$.  If $A$ and $B$ are morphisms in the tangle category with
$Range(B) = Domain(A)$, then the composition of $A$ and $B$ is
denoted $AB$. 

In discussing $Tang$ we shall continue to use the topological
terminology {\em tangle} for a morphism in the category.  An
{\em $n-m$ tangle} is a tangle with
$n$ inputs and $m$ outputs. Thus a $1-1$ tangle is any map from
$V$ to $V$ in
$Tang$, and a $0-0$ tangle is any knot or link arranged with
respect to the vertical to give a morphism from $k$ to $k$. If
$T$ is an $n-m$ tangle in
$Tang$ we define $\Delta (T)$ to be the $2n-2m$ tangle obtained
by replacing every strand of $T$ by two parallel copies of that
strand. 

\subsection{The Functor $F\co Tang \longrightarrow Cat(A)$ }

Let $A$ be a quasi-triangular Hopf algebra. Let $\rho \in A
\otimes A$ denote the Yang--Baxter element for $A$ and write
$\rho$ symbolically in the form $\rho =
\sum e \otimes e'$. We wish to define a functor from $Tang$ to
$Cat(A)$. It suffices to define $F$ on the generating morphisms
$R$,$L$,$Cup$ and $Cap$. We define $$F(Cup) = Cup,$$ $$F(Cap) =
Cap,$$ $$F(R) = P\rho = P (\sum e \otimes e')= \sum P (e
\otimes e'),$$ $$F(L) = \rho^{-1} P = \sum (s(e) \otimes e') P =
\sum (e \otimes s^{-1}(e')) P.$$ 
Diagrammatically, it is convenient to to picture $F(R)$ as a
flat crossing with beads {\em above} the crossing labeled $e$
and $e'$ from left to right, with the summation indicated by
the double appearance of the letter $e$. Similarly,
$F(L)$ is depicted as a crossing with beads {\em below} the
crossing and labeled
$s(e)$ and $e'$.  See Figure 11.

\begin{figure}[ht!]\vglue-2mm
\cl{\small    \setlength{\unitlength}{0.92pt}
\begin{picture}(346,207) \thinlines
\put(270,10){\makebox(25,23){$F(L)$}}
\put(272,108){\makebox(23,21){$F(R)$}}
\put(45,15){\makebox(21,21){$L$}}
\put(46,110){\makebox(20,21){$R$}}
\put(307,36){\makebox(25,25){$e'$}}
\put(229,38){\makebox(25,25){$s(e)$}}
\put(311,157){\makebox(25,25){$e'$}}
\put(234,157){\makebox(22,19){$e$}}
\put(147,68){\makebox(27,21){$F$}}
\put(146,164){\makebox(27,21){$F$}}
\put(115,66){\vector(1,0){96}} \put(114,163){\vector(1,0){96}}
\put(261,48){\circle*{10}} \put(300,48){\circle*{10}}
\put(262,172){\circle*{10}}
\put(303,171){\circle*{10}} \put(237,95){\line(4,-3){85}}
\put(237,31){\line(4,3){85}} \put(240,125){\line(4,3){85}}
\put(240,189){\line(4,-3){85}} \thicklines
\put(11,195){\line(4,-3){86}}
\put(49,157){\line(-4,-3){36}} \put(59,165){\line(4,3){38}}
\put(17,33){\line(4,3){85}} \put(52,67){\line(-4,3){38}}
\put(63,58){\line(4,-3){36}} \end{picture}}\vglue -5mm
\caption{The Functor $F$ on Right and Left
Crossings}
\end{figure}
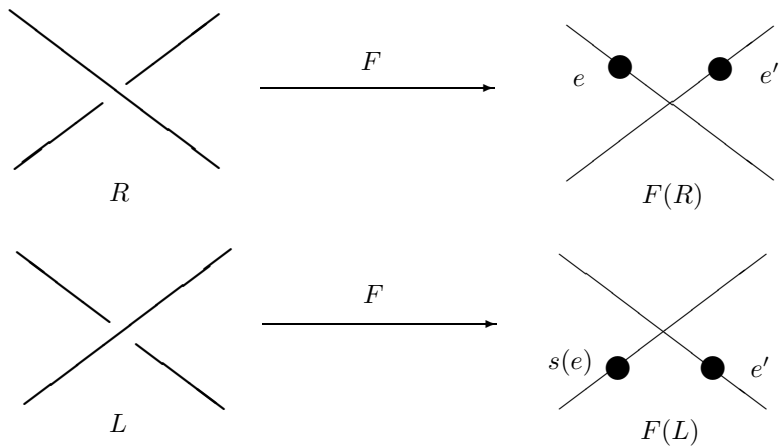

Proof that this functor is well-defined is found in
\cite{KandR} and
\cite{KRS}.

Knowing that $F$ is a functor on the tangle category, means
that we have implicitly defined many invariants of knots and
links, since regularly isotopic tangles will have equivalent
images under $F$.  We can immediately state the following
Theorem for $1-1$ tangles, giving invariants with values in the
Hopf algebra $A$. 

\thm Let $T$ be a single-stranded, $1-1$
tangle. That is , $T$ is a ``knot on a string".  Then $F(T) =
w[T]G^{d(T)}$ is a regular isotopy invariant of $T$.  Here
$w[T]=w(F(T))$ is the element of $A$ defined in the last
section by concentrating all the algebra in $F(T)$ in the lower
part of the tangle, and $d(T)$ is the Whitney degree of the
plane curve underlying $T$.  In fact $w[T] \in A$  is itself a
regular isotopy invariant of $T$, as is $d(T)$.\rm\medskip

\noindent For the proof see \cite{KRS}.\medskip 

{\bf Example}\qua In Figure 12 we illustrate how a
topological curl gives rise under the functor $F$ to the
special central element $v=G^{-1}u$ of a ribbon Hopf algebra.
This element is of particular use in defining the invariants
described in the next section.

\begin{figure}[ht!]\vglue -2mm
\cl{\small    \setlength{\unitlength}{0.92pt}
\begin{picture}(364,292) \thinlines
\put(277,59){\makebox(77,33){$G^{-1}u = v$}}
\put(247,64){\makebox(28,22){$=$}}
\put(123,65){\makebox(28,22){$=$}}
\put(11,64){\makebox(28,22){$=$}}
\put(100,218){\makebox(26,21){$F$}}
\put(79,216){\vector(1,0){76}}
\put(236,194){\makebox(52,25){$s^{-1}(e')$}}
\put(173,199){\makebox(15,16){$e$}}
\put(217,206){\circle*{10}} \thicklines
\put(36,235){\line(1,-1){14}}
\put(70,282){\line(0,-1){46}} \put(70,200){\line(0,-1){46}}
\put(13,199){\line(1,0){21}} \put(12,235){\line(0,-1){35}}
\put(35,235){\line(-1,0){21}} \put(35,200){\line(1,1){34}}
\put(55,215){\line(1,-1){15}} \thinlines
\put(223,282){\line(0,-1){47}}
\put(223,234){\line(-1,-1){36}} \put(186,197){\line(-1,0){21}}
\put(165,197){\line(0,1){37}} \put(166,234){\line(1,0){22}}
\put(189,234){\line(1,-1){35}} \put(224,198){\line(0,-1){46}}
\put(109,59){\line(0,-1){46}} \put(74,95){\line(1,-1){35}}
\put(51,95){\line(1,0){22}} \put(50,58){\line(0,1){37}}
\put(71,58){\line(-1,0){21}} \put(108,95){\line(-1,-1){36}}
\put(108,143){\line(0,-1){47}} \put(226,56){\line(0,-1){46}}
\put(191,92){\line(1,-1){35}} \put(168,92){\line(1,0){22}}
\put(167,55){\line(0,1){37}} \put(188,55){\line(-1,0){21}}
\put(225,92){\line(-1,-1){36}} \put(225,140){\line(0,-1){47}}
\put(51,76){\circle*{10}} \put(226,108){\circle*{10}}
\put(225,126){\circle*{10}}
\put(109,124){\circle*{10}} \put(196,206){\circle*{10}}
\put(115,116){\makebox(15,16){$e$}}
\put(232,118){\makebox(15,16){$e$}}
\put(57,67){\makebox(20,19){$e'$}}
\put(234,96){\makebox(40,23){$s(e'$)}}
\end{picture}}\vglue -2mm
\caption{The Ribbon Element $v$}
\end{figure}
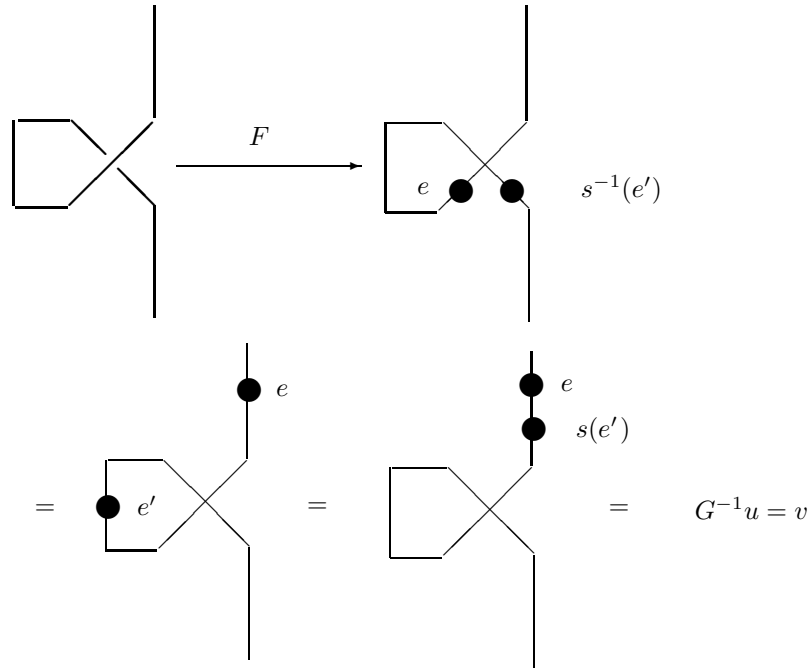

\section {Invariants of 3--manifolds}

The structure we have built so far can be used to construct
invariants of 3--manifolds presented in terms of surgery on
framed links.  We sketch here our technique that simplifies  an
approach to 3--manifold invariants of Mark Hennings
\cite{Hennings}. 

Recall that an element  $\lambda$  of  the dual algebra
$A^{\ast}$  is said to be a {\em right integral }  if $\lambda
(x)1 = m(\lambda \otimes 1)(\Delta (x))$ for all $x$  in $A.$
For a unimodular  \cite{LarsonSweedler} ,\cite{Radford-trace}
finite dimensional  ribbon Hopf algebra  $A$ there is a  right
integral
$\lambda$ satisfying the following properties  for all x and y
in A:

\noindent 0)\qua $\lambda$ is unique up to scalar multiplication
when $k$  is a field. 

\noindent 1)\qua $\lambda (xy) = \lambda (s^{2}(y)x).$ 

\noindent 2)\qua $\lambda (gx) = \lambda (s(x))$  where $g=
G^{2}$,  $G$ the special grouplike element for the ribbon
element  $v= G^{-1}u.$

Given the existence of this  integral $\lambda$,  define a
functional  $tr\co A
\longrightarrow k$ by the formula $$tr(x) =  \lambda (Gx).$$
(It follows from the fact that $s^{2}(G)=G$ that $\lambda(Gx) =
\lambda(xG).$) 

It is then easy to prove the following theorem \cite{KandR}.

\begin{tthm} The function $tr$ defined as
above satisfies
$$tr(xy) = tr(yx)\hbox{ for all } x,y\hbox{ in }A$$ 
$$tr(s(x)) = tr(x)\hbox{ for all } x,y\hbox{ in }A.\leqno{and}$$
\end{tthm}

The upshot of this theorem is that for a unimodular finite
dimensional Hopf algebra there is a natural trace defined via
the existent right integral. Remarkably, this  trace is just
designed  to behave well with respect to handle sliding
\cite{K-Hopf}, \cite{KandR}.  Handle sliding  is the basic
transformation on framed links that leaves the corresponding
3--manifold obtained by framed surgery unchanged.  See
\cite{KirbyCalc}.  This means that a suitably normalized
version of this trace on framed links gives an invariant of
3--manifolds. For a link $K$, we let TR(K) denote the functional
on links, as described in  the previous section, defined via
$tr$ as above. 

To see how the condition on handle sliding and the property of
being a right integral are related in our category, we refer
the reader to Figure 13 where the basic form of handle sliding
is illustrated and its algebraic counterpart is shown. The
algebraic counterpart arises when we concentrate all the algebra
in a given link component in one place on the diagram. The
component is then replaced by a circle and formally its
evaluation is $O_{R}(x)$ for a suitable $x$ in the Hopf algebra
where $O_{R}(x)$ denotes the evaluation of the morphism
corresponding to a circle labelled with $x$ as in the Figure.
As the diagram shows, if we let $\lambda(x) = O_{R}(x)$, then
invariance under handle sliding is implicated by $\lambda$
being a right integral on the Hopf algebra.

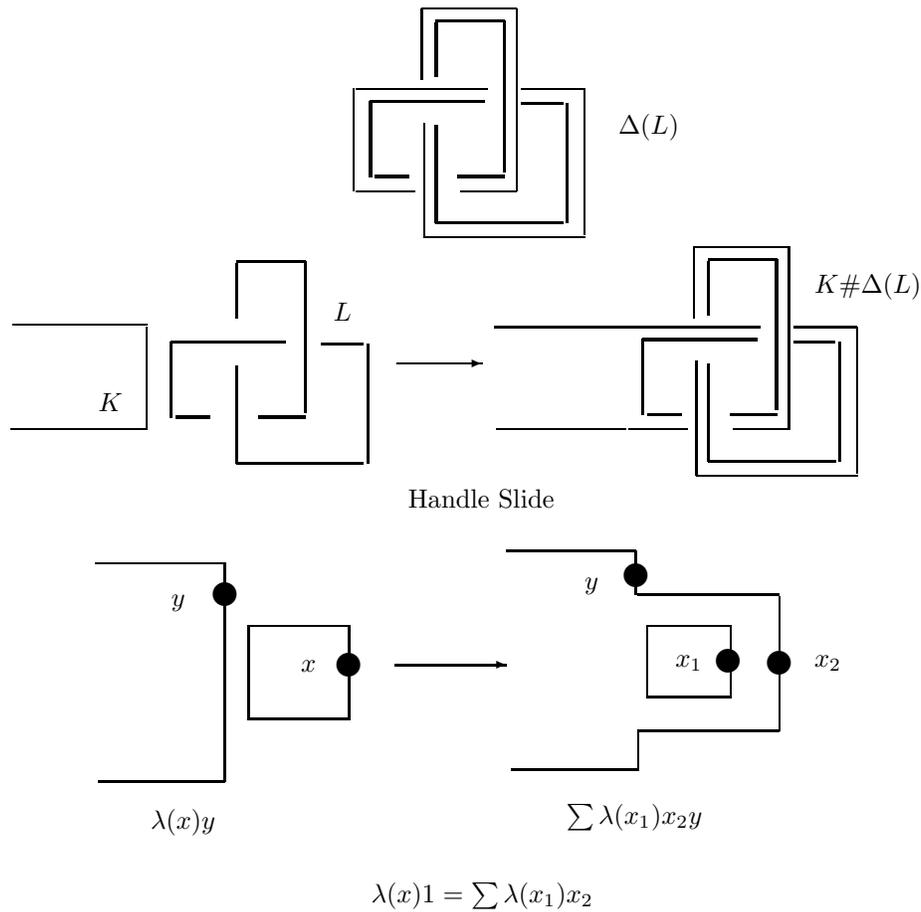
\begin{figure}[ht!]\vglue -2mm
\cl{\small    \setlength{\unitlength}{0.92pt}
\begin{picture}(401,409) \thinlines
\put(115,10){\makebox(178,44){$\lambda(x)1 = \sum
\lambda(x_{1})x_{2}$}}
\put(98,156){\circle*{10}} \thicklines
\put(76,260){\line(1,0){47}}
\put(138,259){\line(1,0){17}} \put(155,210){\line(-1,0){52}}
\put(103,210){\line(0,1){40}} \put(103,270){\line(0,1){22}}
\put(103,293){\line(1,0){28}} \put(131,293){\line(0,-1){62}}
\put(131,229){\line(-1,0){19}} \put(92,229){\line(-1,0){14}}
\put(76,260){\line(0,-1){31}} \put(157,259){\line(0,-1){49}}
\thinlines
\put(264,266){\line(1,0){54}} \put(66,224){\line(0,1){42}}
\put(288,224){\line(-1,0){24}} \put(330,224){\line(-1,0){23}}
\put(330,299){\line(0,-1){75}} \put(291,299){\line(1,0){39}}
\put(291,270){\line(0,1){29}} \put(292,205){\line(0,1){47}}
\put(358,205){\line(-1,0){66}} \put(358,266){\line(0,-1){60}}
\put(332,266){\line(1,0){26}} \thicklines
\put(270,261){\line(1,0){47}}
\put(332,260){\line(1,0){17}} \put(349,211){\line(-1,0){52}}
\put(297,211){\line(0,1){40}} \put(297,271){\line(0,1){22}}
\put(297,294){\line(1,0){28}} \put(325,294){\line(0,-1){62}}
\put(325,230){\line(-1,0){19}} \put(286,230){\line(-1,0){14}}
\put(270,261){\line(0,-1){31}} \put(351,260){\line(0,-1){49}}
\thinlines
\put(264,266){\line(-1,0){55}} \put(263,224){\line(-1,0){53}}
\put(66,267){\line(-1,0){55}} \put(66,224){\line(-1,0){56}}
\put(268,156){\line(1,0){58}} \put(326,155){\line(0,-1){55}}
\put(326,100){\line(-1,0){58}} \put(214,174){\line(1,0){53}}
\put(267,174){\line(0,-1){18}} \put(216,84){\line(1,0){52}}
\put(268,85){\line(0,1){14}} \put(46,79){\line(1,0){52}}
\put(45,169){\line(1,0){53}} \put(108,105){\framebox(41,38){}}
\put(98,169){\line(0,-1){90}} \put(267,164){\circle*{10}}
\put(149,127){\circle*{10}} \put(326,128){\circle*{10}}
\put(305,129){\circle*{10}} \put(272,114){\framebox(34,29){}}
\put(169,251){\vector(1,0){35}}
\put(148,183){\makebox(111,25){Handle Slide}}
\put(41,227){\makebox(20,16){$K$}}
\put(136,263){\makebox(21,19){$L$}}
\put(334,270){\makebox(57,26){$K\#\Delta(L)$}}
\put(168,127){\vector(1,0){46}}
\put(124,118){\makebox(17,17){$x$}}
\put(68,141){\makebox(22,23){$y$}}
\put(238,148){\makebox(22,23){$y$}}
\put(280,120){\makebox(18,16){$x_{1}$}}
\put(333,117){\makebox(26,22){$x_{2}$}}
\put(251,333){\makebox(43,31){$\Delta(L)$}}
\put(152,364){\line(1,0){54}}
\put(151,322){\line(0,1){42}} \put(176,322){\line(-1,0){24}}
\put(218,322){\line(-1,0){23}} \put(218,397){\line(0,-1){75}}
\put(179,397){\line(1,0){39}} \put(179,368){\line(0,1){29}}
\put(180,303){\line(0,1){47}} \put(246,303){\line(-1,0){66}}
\put(246,364){\line(0,-1){60}} \put(220,364){\line(1,0){26}}
\thicklines
\put(158,359){\line(1,0){47}} \put(220,358){\line(1,0){17}}
\put(237,309){\line(-1,0){52}} \put(185,309){\line(0,1){40}}
\put(185,369){\line(0,1){22}} \put(185,392){\line(1,0){28}}
\put(213,392){\line(0,-1){62}} \put(213,328){\line(-1,0){19}}
\put(174,328){\line(-1,0){14}} \put(158,359){\line(0,-1){31}}
\put(239,358){\line(0,-1){49}} \thinlines
\put(28,44){\makebox(106,37){$\lambda(x)y$}}
\put(202,44){\makebox(129,41){$\sum\lambda(x_{1})x_{2}y$}}
\end{picture}}\vglue-2mm
\caption{Handle Sliding and Right Integral}
\end{figure}

A proper normalization of $TR(K)$ gives an invariant of the
3--manifold obtained by framed surgery on $K.$  More precisely
(assuming that $\lambda(v)$ and
$\lambda(v^{-1})$ are non-zero),  let
$$INV(K) = (\lambda(v) \lambda(v^{-1}))^{- c(K)/2}(\lambda(v) /
\lambda(v^{-1})^{- \sigma(K)/2}TR(K)$$
\noindent where  $c(K)$ denotes the number of components of K,
and
$\sigma(K)$ denotes the signature of the matrix of linking
numbers of the components of $K$ (with framing numbers on the
diagonal).  Then  $INV(K)$  is an invariant of the 3--manifold
obtained by doing framed surgery on $K$ in the blackboard
framing. This is our reconstruction  of Hennings invariant
\cite{Hennings} in an intrinsically unoriented context.

\Addresses\recd

\end{document}